\documentclass[p1]{elsarticle}
\usepackage[utf8]{inputenc}
\usepackage{amsfonts}
\usepackage{amssymb}
\usepackage{amsmath}
\usepackage{hyperref}
\hypersetup{colorlinks = true, citecolor = blue, linkcolor = blue}

\title{A sum up method for solving summations of the form $\sum_{k=n_0}^{n} A_{n,k}$ and rising and falling factorial transforms}
\author{Parham Zarghami\corref{cor1}\fnref{fn1}} 
\cortext[cor1]{Corresponding author \ead{parham.zarghami@ut.ac.ir}}
\fntext[fn1]{Declarations of interest: None}

\affiliation{organization={University of Tehran}, addressline={School of Electrical and Computer Engineering, University College of Engineering, North Kargar street}, city={Tehran}, postcode={1439957131}, country={Iran}}

\date{December 2023}

\begin{document}

\begin{abstract}
    In this paper, we discuss a method that utilizes the recurrence of $A_{n,k}$ to solve summations of the form $\sum_{k=n_0}^{n} A_{n,k}$. It is observed that by repeating the procedure, the upper bound of summation is reduced and tilts toward the lower bound. This method of summation is mostly suitable for combinatorial sequences such as binomial coefficients, Stirling numbers of both kinds, etc. After the main method is displayed, some examples are illustrated. Some useful identities about Stirling and r-Stirling numbers are obtained. Finally, two transforms called rising and falling factorial transforms which turn the basis of power polynomials into factorial basis are derived. These transforms verify and simplify the results obtained in the examples section. Also, these transforms describe the relationship between fractional derivatives (or fractional integrals) and falling factorial (or rising factorial) by its series expansion.
\end{abstract}


\maketitle
\section{Introduction}
There are several formulas and methods, e.g., telescoping series, that focus on solving special summations. A summation is a way of accumulating information in an explicit form. However, some information may be similar and can be summed up before the summation. The aim of this research is to find a method of summing up the same terms utilizing recurrences between different terms of the sequence. This way, the bound of summations dwindles away until it becomes trivial.

Table \ref{tab:Table 1} shows a number of sequences with their recurrence relations. Table \ref{tab:Table 2} illustrates definition of required functions.

\begin{table}[h!]
    \caption{Recurrence relation of sequences}
    \label{tab:Table 1}
    \centering
    \small
    \begin{tabular}{|c|c|}
      \hline
      \textbf{Sequence} & \textbf{Recurrent relation}\\
      \hline
      Binomial coefficients (A007318) & $\binom{n}{k}=\binom{n-1}{k}+\binom{n-1}{k-1}$\\
      \hline
      Stirling numbers of the first kind (A008275) & ${n\brack k} = (n-1){n-1\brack k} + {n-1\brack k-1}$\\
      \hline
      Stirling numbers of the second kind (A008277) & ${n\brace k} = k{n-1\brace k} + {n-1\brace k-1}$\\
      \hline
      r-Stirling numbers of the first kind \cite{ref. 7} & ${n\brack k}_r = (n-1){n-1\brack k}_r + {n-1\brack k-1}_r$\\
      \hline
      r-Stirling numbers of the second kind \cite{ref. 7} & ${n\brace k}_r = k{n-1\brace k}_r + {n-1\brace k-1}_r$\\
      \hline
    \end{tabular}
\end{table}

\begin{table}[h!]
    \caption{Definition of required functions}
    \label{tab:Table 2}
    \centering
    \small
    \begin{tabular}{|c|c|}
        \hline
        \textbf{Function} & \textbf{Definition}\\
        \hline
        Rising factorial \cite{ref. 1,ref. 2} & $x^{\overline{n}}= \prod_{i=0}^{n-1} (x+i)=\frac{\Gamma (x+n)}{\Gamma (x)}$\\
        \hline
        Falling factorial \cite{ref. 1,ref. 2} & $(x)_n= \prod_{i=0}^{n-1} (x-i)=\frac{\Gamma (x+1)}{\Gamma (x-n+1)}$\\
        \hline
        Touchard (or Bell) polynomials \cite{ref. 9} & $T_n(x) = \sum_{k=0}^{n} {n\brace k} x^k $\\
        \hline
        r-Touchard polynomials \cite{ref. 9} & $T_{n,r}(x) = \sum_{k=0}^{n} {n+r\brace k+r}_r x^k $\\
        \hline
        Incomplete gamma function \cite{ref. 1} & $\Gamma(s,x) = \int_{x}^{\infty} t^{s-1}e^{-t}dt$\\
        \hline
    \end{tabular}
\end{table}

The following formulas show some features of incomplete gamma function \cite{ref. 1}:
\begin{equation}
\label{eq. 1}
    \Gamma(n+1,x) = \int_{x}^{\infty} t^{n}e^{-t}dt = \sum_{k=0}^{n} \frac{n!}{k!} x^k e^{-x}
\end{equation}
\begin{equation}
    \Gamma(s+1,x) = s\Gamma(s+1,x) + x^s e^{-x}
\end{equation}

In section 2, we will propose our method and show how some specific summations tend to assume an explicit form through recurrences. In section 3, some applications of our proposed method are demonstrated. In section 4, two new transformation named rising and falling factorial transforms that can change the basis of polynomial and series from power to factorial ones is proposed. Finally, suggestions for future works and conclusions are presented.

\section{The method}

Our proposed method is similar to the Telescoping series. The difference is in summing up the terms instead of eliminating the similar terms. To illustrate how our method works, we consider an example recurrence:
\begin{equation}
    A_{n,k}=\sum_{i=0}^{m} a_i(n,k)A_{n-1,k-i}
\end{equation}
Summation on $k$:
\begin{equation}
    \sum_{k=n_0}^{n} A_{n,k}=\sum_{k=n_0}^{n} \sum_{i=0}^{m} a_i(n,k)A_{n-1,k-i} = \sum_{k=n_0}^{n-1} c_{1,i}(n)A_{n-1,i}
\end{equation}
Then,  $A_{n,m}=0$ and $a_i(n,m)=0$ for $n<m$ and $m<n_0$ is considered. Based on above process, $c_{1,i}(n)$ can be defined as follows:
\begin{equation}
    c_{1,i}(n)=\sum_{j=0}^{k} a_j(n,i+j)
\end{equation}
This procedure is repeated as the following:
\begin{equation}
    c_{1,k}(n)A_{n,k}=\sum_{i=0}^{m} a_i(n,k)A_{n-1,k-i}
\end{equation}
\begin{equation}
    \sum_{k=n_0}^{n} A_{n,k} = \sum_{k=n_0}^{n-1} c_{1,k}(n)A_{n-1,k} = \sum_{k=n_0}^{n-2} c_{2,k}(n)A_{n-2,k}
\end{equation}
Again based on above processes, $c_{2,i}(n)$ can be formulated as:
\begin{equation}
    c_{2,i}(n)=\sum_{j=0}^{k} c_{1,i+j}(n)a_j(n-1,i+j)
\end{equation}
After that, the following formulations can be constructed:
\begin{equation}
\label{eq. 9}
    \sum_{k=n_0}^{n} A_{n,k} = \sum_{k=n_0}^{n-1} c_{1,k}(n)A_{n-1,k} = \sum_{k=n_0}^{n-2} c_{2,k}(n)A_{n-2,k}= ... = c_{n-n_0,n_0}(n)A_{n_0,n_0}
\end{equation}
\begin{equation}
\label{eq. 10}
    c_{s,i}(n)=\sum_{j=0}^{k} c_{s-1,i+j}(n)a_j(n-s+1,i+j)
\end{equation}
This method may be applied to the super-recurrence of the form:
\begin{equation}
    \left|\genfrac{}{}{0pt}{}{n}{k}\right|=f(n,k)\left|\genfrac{}{}{0pt}{}{n-1}{k}\right|+g(n,k)\left|\genfrac{}{}{0pt}{}{n-1}{k-1}\right|+[n=k=0]
\end{equation}
For instance in \cite{ref. 3}, the coefficients used are $f(n,k)=\alpha n+\beta k+\gamma , g(n,k)=\alpha^{'} n+\beta^{'} k+\gamma^{'}$. By the definition, $\left|\genfrac{}{}{0pt}{}{0}{0}\right|=1$ and $\left|\genfrac{}{}{0pt}{}{n}{k}\right|=0$ when $n<k$. So, the value of $\left|\genfrac{}{}{0pt}{}{n}{n}\right|$ in the following form can be calculated by the above recurrence:
\begin{equation}
	\begin{split}
		\left|\genfrac{}{}{0pt}{}{n}{n}\right| & = g(n,n)\left|\genfrac{}{}{0pt}{}{n-1}{n-1}\right| \\
		& =g(n,n)g(n-1,n-1)\left|\genfrac{}{}{0pt}{}{n-2}{n-2}\right|= ... =\prod_{i=1}^{n} g(i,i)
	\end{split} 
\end{equation}
(\ref{eq. 9}) and (\ref{eq. 10}) are used for polynomials of the form $\sum_{i=n_0}^{n} \left|\genfrac{}{}{0pt}{}{n}{k}\right| x^k $ as follows:
\begin{equation}
    \sum_{i=n_0}^{n} \left|\genfrac{}{}{0pt}{}{n}{k}\right| x^k = x^{n_0} \left|\genfrac{}{}{0pt}{}{n_0}{n_0}\right| Y_{n-n_0,n_0}(x) = x^{n_0} Y_{n-n_0,n_0}(x) \prod_{i=1}^{n_0} g(i,i)
\end{equation}
Where $Y_{m,k}(x)$ has a recurrence as below:
\begin{equation}
    Y_{m,k}(x)=f(n-m+1,k)Y_{m-1,k}(x)+xg(n-m+1,k+1)Y_{m-1,k+1}(x)
\end{equation}
Here $n$ is the upper bound of the summation and $Y_{0,k}(x)=1$, $Y_{1,k}(x)=f(n,k)+x g(n,k+1)$.
(\ref{eq. 9}) and (\ref{eq. 10}) are used for rising factorial polynomials of the form $\sum_{k=n_0}^{n} \left|\genfrac{}{}{0pt}{}{n}{k}\right| x^{\overline{k}}$ as follows:
\begin{equation}
    \sum_{k=n_0}^{n} \left|\genfrac{}{}{0pt}{}{n}{k}\right| x^{\overline{k}} = x^{\overline{n_0}} \left|\genfrac{}{}{0pt}{}{n_0}{n_0}\right| y_{n-n_0,n_0}(x) = x^{\overline{n_0}} y_{n-n_0,n_0}(x) \prod_{i=1}^{n_0} g(i,i)
\end{equation}
Where $y_{m,k}(x)$ has a recurrence as below:
\begin{equation}
\label{eq. 16}
    y_{m,k}(x)=f(n-m+1,k)y_{m-1,k}(x) + (x+k)g(n-m+1,k+1)y_{m-1,k+1}(x)
\end{equation}
$g(n,n_0)\left|\genfrac{}{}{0pt}{}{n-1}{n_0 - 1}\right|=0$ is considered at the middle of calculations of both polynomials. In (\ref{eq. 16}), $n$ is the upper bound of the summation and $y_{0,k}(x)=1$, $y_{1,k}(x)=f(n,k)+(x+k)g(n,k+1)$.

\section{Examples}

In this section, we will focus on special cases of super-recurrence involved in polynomials discussed in the previous section. Table \ref{tab:Table 3} illustrates recurrence of three famous combinatorial sequences:
\begin{table}[h!]
    \caption{Recurrences of sequences}
    \label{tab:Table 3}
    \centering
    \small
    \begin{tabular}{|c|c|c|c|}
      \hline
      \textbf{Sequence} & \textbf{Recurrence} & \textbf{$f(n,k)$} & \textbf{$g(n,k)$}\\
      \hline
      Binomial coefficients & $\left|\genfrac{}{}{0pt}{}{n}{k}\right|=\left|\genfrac{}{}{0pt}{}{n-1}{k}\right|+\left|\genfrac{}{}{0pt}{}{n-1}{k-1}\right|$ & 1 & 1\\
      \hline
      Stirling numbers of the first kind & $\left|\genfrac{}{}{0pt}{}{n}{k}\right|=(n-1) \left|\genfrac{}{}{0pt}{}{n-1}{k}\right|+\left|\genfrac{}{}{0pt}{}{n-1}{k-1}\right|$ & n-1 & 1\\
      \hline
      Stirling numbers of the second kind & $\left|\genfrac{}{}{0pt}{}{n}{k}\right|=k \left|\genfrac{}{}{0pt}{}{n-1}{k}\right|+\left|\genfrac{}{}{0pt}{}{n-1}{k-1}\right|$ & k & 1\\
      \hline
    \end{tabular}
\end{table}

Recurrence relation of $Y_{m,k}(x)$ and $y_{m,k}(x)$ illustrated as Table \ref{tab:Table 4} and Table \ref{tab:Table 5} respectively.

\begin{table}[h!]
    \caption{$Y_{m,k}(x)$ of sequences}
    \label{tab:Table 4}
    \centering
    \small
    \begin{tabular}{|c|c|}
      \hline
      \textbf{Sequence} & \textbf{$Y_{m,k}(x)$}\\
      \hline
      Binomial coefficients & $Y_{m,k}(x)=Y_{m-1,k}(x)+x Y_{m-1,k+1}(x)$\\
      \hline
      Stirling numbers of the first kind & $Y_{m,k}(x)=(n-m) Y_{m-1,k}(x)+x Y_{m-1,k+1}(x)$\\
      \hline
      Stirling numbers of the second kind & $Y_{m,k}(x)=k Y_{m-1,k}(x)+x Y_{m-1,k+1}(x)$\\
      \hline
    \end{tabular}
\end{table}

\begin{table}[h!]
    \caption{$y_{m,k}(x)$ of sequences}
    \label{tab:Table 5}
    \centering
    \small
    \begin{tabular}{|c|c|}
      \hline
      \textbf{Sequence} & \textbf{$y_{m,k}(x)$}\\
      \hline
      Binomial coefficients & $y_{m,k}(x)=y_{m-1,k}(x)+(x+k) y_{m-1,k+1}(x)$\\
      \hline
      Stirling numbers of the first kind & $y_{m,k}(x)=(n-m) y_{m-1,k}(x)+(x+k) y_{m-1,k+1}(x)$\\
      \hline
      Stirling numbers of the second kind & $y_{m,k}(x)=k y_{m-1,k}(x)+(x+k) y_{m-1,k+1}(x)$\\
      \hline
    \end{tabular}
\end{table}

It is observed that in all cases $Y_{0,k}(x)=1$ and $y_{0,k}(x)=1$ are present. For solving the above recurrences, we do not need to use complicated methods while we can use the recurrence unfolding method \cite{ref. 3}. It means that we must replace $n$th term with $(n-1)$th terms involved in the main recurrence. We can perform iterations in a top-down and bottom-up manner. To avoid repeated calculations, we have displayed the results only.

\subsection{Solving $Y_{m,k}(x)=Y_{m-1,k}(x)+x Y_{m-1,k+1}(x)$}

$Y_{m,k}(x)$ is independent of $k$, so the recurrence takes the form $Y_{m,k}(x)=(x+1) Y_{m-1,k}(x)$ form. Hence:
\begin{equation}
    Y_{n,k}(x)= (x+1)^n
\end{equation}
So, for the binomial coefficients involved in following polynomial:
\begin{equation}
    \sum_{i=0}^{n} \binom{n}{i} x^i = Y_{n,0}(x) = (x+1)^n
\end{equation}

\subsection{Solving $y_{m,k}(x)=y_{m-1,k}(x)+(x+k) y_{m-1,k+1}(x)$}

First, we want to show the results by using top-down unfolding method:
\begin{equation}
    y_{n,k}(x) = \sum_{i=0}^{n-1} \binom{n-1}{i} (x+k+i+1) (x+k)^{\overline{i}}
\end{equation}
We know that $\Delta_x ((x+k)^{\overline{i}}) = i(x+k+1)^{\overline{i-1}} = i \frac{(x+k)^{\overline{i}}}{x+k}$:
\begin{equation}
    y_{n,k}(x)=(x+k+1)y_{n-1,0}(x+k) + (x+k)\Delta_x y_{n-1,0}(x+k+1)
\end{equation}
Is simplified to:
\begin{equation}
\label{eq. 21}
    y_{n,k}(x)=y_{n-1,0}(x+k) + (x+k) y_{n-1,0}(x+k+1)=y_{n,0}(x+k)
\end{equation}
We used $y_{n,0}(x)=\sum_{i=0}^{n} \binom{n}{i} x^{\overline{i}}$ to obtain (\ref{eq. 21}). By using $x=1$ and (\ref{eq. 1}):
\begin{equation}
	\begin{split}
		y_{n,0}(1) & =\sum_{i=0}^{n} \binom{n}{i} i!=\sum_{i=0}^{n} \frac{n!}{i!} \\
		& =e \Gamma (n+1,1)=e\int_{1}^{\infty} t^n e^{-t} dt
	\end{split}
\end{equation}
Without loss of generality, in (\ref{eq. 21}), we can consider $k=0$ and $x=m$ which $m$ is an integer:
\begin{equation}
    my_{n,0}(m+1)=y_{n+1,0}(m)-y_{n,0}(m) = \Delta_n y_{n,0}(m)
\end{equation}
Taking $m$ difference on $n$ (i.e. $\Delta_{n}^{m-1}$) from both side of equation results in:
\begin{equation}
    y_{n,0}(m)=\frac{1}{(m-1)!} \Delta_{n}^{m-1} y_{n,0}(1) = \frac{e}{(m-1)!} \int_{1}^{\infty} \Delta_{n}^{m-1} (t^n) e^{-t} dt
\end{equation}
As the difference of exponential function is $\Delta_{n}^{m-1} t^n=(t-1)^{m-1} t^n$:
\begin{equation}
    y_{n,0}(m)=\frac{e}{(m-1)!} \int_{1}^{\infty} (t-1)^{m-1} t^n e^{-t} dt
\end{equation}
$m$ is an integer variable. Thus, a general form like the following is needed:
\begin{equation}
    y_{n,0}(x)=\frac{1}{\Gamma (x)} \int_{0}^{\infty} (t+1)^{n} t^{x-1} e^{-t} dt
\end{equation}
This verifies main recurrence. Then we arrive in an explicit form for $y_{n,k}(x)$ by (\ref{eq. 21}):
\begin{equation}
    y_{n,k}(x)=y_{n,0}(x+k) = \frac{1}{\Gamma (x+k)} \int_{0}^{\infty} (t+1)^{n} t^{x+k-1} e^{-t} dt
\end{equation}
And integral representation for $\sum_{k=0}^{n} \binom{n}{k} x^{\overline{k}}$:
\begin{equation}
	\begin{split}
		\sum_{k=0}^{n} \binom{n}{k} x^{\overline{k}} & = y_{n,0}(m) \\
		& = \frac{1}{\Gamma (x)} \int_{0}^{\infty} (t+1)^{n} t^{x-1} e^{-t} dt
	\end{split}
\end{equation}

\subsection{Solving $Y_{m,k}(x)=(n-m) Y_{m-1,k}(x)+x Y_{m-1,k+1}(x)$}

With unfolding method we can find the fact that $Y_{m,k}(x)$ is independent from $k$ similar to the recurrence that we discussed in 3.1. So, this recurrence takes the following form:
\begin{equation}
    Y_{m,k}(x)=(x+n-m)Y_{m-1,k}(x)
\end{equation}
We conclude that:
\begin{equation}
    Y_{m,k}(x)=(x+n-m)^{\overline{m}}
\end{equation}
Then, the explicit form of intended polynomial is easily obtained:
\begin{equation}
    \sum_{k=0}^{n} {n \brack k} x^k = Y_{n,0}(x) = x^{\overline{n}}
\end{equation}

\subsection{Solving $y_{m,k}(x)=(n-m) y_{m-1,k}(x)+(x+k) y_{m-1,k+1}(x)$}

First of all, we should change $y_{m,k}(x)$ to $s_{m,k}(x)$ as below:
\begin{equation}
\label{eq. 32}
    s_{m,k}(x)=\Gamma (x+k) y_{m,k}(x)  , s_{0,k}(x) = \Gamma (x+k)
\end{equation}
Next, using bottom-up unfolding method we get:
\begin{equation}
    s_{m,k}(x) = \sum_{j=0}^{m} \sum_{i=0}^{j} {m \brack j} \binom{j}{i} (-1)^{m-j} (n-1)^{j-i} \Gamma (x+k+i)
\end{equation}
Applying definition of gamma function:
\begin{equation}
    s_{m,k}(x) = \sum_{j=0}^{m} \sum_{i=0}^{j} {m \brack j} \binom{j}{i} (-1)^{m-j} (n-1)^{j-i} \int_{0}^{\infty} t^{x+k+i-1} e^{-t} dt
\end{equation}
It is simplified as:
\begin{equation}
\begin{split}
    s_{m,k}(x) & = \int_{0}^{\infty} \sum_{j=0}^{m} (\sum_{i=0}^{j} {m \brack j} \binom{j}{i} (-1)^{m-j} (n-1)^{j-i} t^{x+k+i-1}) e^{-t} dt \\
    & = \int_{0}^{\infty} \sum_{j=0}^{m} {m \brack j} (-1)^{m-j} (t+n-1)^{j} t^{x+k-1} e^{-t} dt \\
    & = \int_{0}^{\infty} (t+n-1)_m t^{x+k-1} e^{-t} dt
\end{split}
\end{equation}
And by (\ref{eq. 32}) we get the intended explicit formula:
\begin{equation}
    y_{m,k}(x) = \frac{1}{\Gamma (x+k)} \int_{0}^{\infty} (t+n-1)_m t^{x+k-1} e^{-t} dt
\end{equation}
Suppose $k=0$ and $m=n$:
\begin{equation}
	\begin{split}
		y_{n,0}(x) & = \frac{1}{\Gamma (x)} \int_{0}^{\infty} (t+n-1)_n t^{x-1} e^{-t} dt \\
		& = \frac{1}{\Gamma (x)} \int_{0}^{\infty} t^{\overline{n}} t^{x-1} e^{-t} dt
	\end{split}
\end{equation}
Then for the Stirling number of the first kind involved in rising factorial polynomial we have the formula below:
\begin{equation}
    \sum_{k=0}^{n} {n \brack k} x^{\overline{k}} = y_{n,0}(x) = \frac{1}{\Gamma (x)} \int_{0}^{\infty} t^{\overline{n}} t^{x-1} e^{-t} dt
\end{equation}

\subsection{Solving $Y_{m,k}(x)=k Y_{m-1,k}(x)+x Y_{m-1,k+1}(x)$}

By top-down unfolding method we can derive the following formula:
\begin{equation}
\label{eq. 39}
    Y_{n,k}(x)= \sum_{i=0}^{n} {n+k \brace i+k}_k x^i = T_{n,k}(x)
\end{equation}
In which $T_{r,k}(x)$ is r-Touchard polynomials shown in section 1. From \cite{ref. 7} we know that:
\begin{equation}
    {n+k \brace i+k}_k = \frac{1}{i!} \sum_{j=0}^{i} \binom{i}{j} (-1)^{i-j} (k+j)^n
\end{equation}
So:
\begin{equation}
\begin{split}
    Y_{n,k}(x) & = \sum_{i=0}^{n} {n+k \brace i+k}_k x^i \\
    & = \sum_{i=0}^{n} \frac{1}{i!} \sum_{j=0}^{i} \binom{i}{j} (-1)^{i-j} (k+j)^n x^i \\
    & = \sum_{j=0}^{n} \frac{(k+j)^n}{j!} \sum_{i=j}^{n} (-1)^{i-j}  \frac{x^i}{(i-j)!} \\
    & = \sum_{j=0}^{n} \frac{(k+j)^n x^j}{j!} \sum_{i=0}^{n-j} \frac{(-x)^i}{i!}
\end{split}
\end{equation}
By (\ref{eq. 1}):
\begin{equation}
\begin{split}
    Y_{n,k}(x) & = \frac{e^{-x}}{n!} \sum_{j=0}^{n} \binom{n}{j} (k+j)^n x^j \int_{-x}^{\infty} t^{n-j} e^{-t} dt \\
    & = \frac{e^{-x}}{n!} \int_{-x}^{\infty} \sum_{j=0}^{n} \binom{n}{j} (k+j)^n x^j t^{n-j} e^{-t} dt
\end{split}
\end{equation}
Suppose $D_x$ is an operator in $D_x=x \frac{d}{dx}$, then, we can obtain $D_x^n (x^i)=i^n x^i$. So:
\begin{equation}
\begin{split}
    Y_{n,k}(x) & = \frac{(-1)^n x^{-k} e^{-x}}{n!} \int_{-x}^{\infty} \sum_{j=0}^{n} \binom{n}{j} (-k-j)^n x^{k+j} t^{-k-j} t^{n+k} e^{-t} dt \\
    & = \frac{(-1)^n x^{-k} e^{-x}}{n!} \int_{-x}^{\infty} \sum_{j=0}^{n} \binom{n}{j} D_{t}^{n} ((\frac{x}{t})^{k+j}) t^{n+k} e^{-t} dt \\
    & = \frac{(-1)^n x^{-k} e^{-x}}{n!} \int_{-x}^{\infty} D_{t}^{n} ((\frac{x}{t})^k (1+\frac{x}{t})^n) t^{n+k} e^{-t} dt
\end{split}
\end{equation}
The operator in integral can be vanished with integration by part. Also, this integral can be solved by Laplace transform or other transforms. For intended polynomial we have:
\begin{equation}
    \sum_{k=0}^{n} {n \brace k} x^k = T_n (x) = Y_{n,0}(x) = \frac{(-1)^n e^{-x}}{n!} \int_{-x}^{\infty} D_{t}^{n} ( (1+\frac{x}{t})^n) t^{n} e^{-t} dt
\end{equation}
The same result can be obtained by Dobinski formula \cite{ref. 5} as follows:
\begin{equation}
    \sum_{k=0}^{n} \frac{k^n}{k!} x^k = \frac{(-1)^n}{n!} \int_{-x}^{\infty} D_{t}^{n} ( (1+\frac{x}{t})^n) t^{n} e^{-t} dt
\end{equation}
Also for Bell numbers \cite{ref. 4}:
\begin{equation}
    B_n = \frac{1}{e} \sum_{k=0}^{n} \frac{k^n}{k!} = \frac{(-1)^n}{e.n!} \int_{-1}^{\infty} D_{t}^{n} ( (1+\frac{1}{t})^n) t^{n} e^{-t} dt
\end{equation}
Using bottom up unfolding for this recurrence on $k$ and applying $Y_{n,0}(x)=T_n (x)$, we can easily observe that:
\begin{equation}
    Y_{n,k}(x)= \frac{1}{x^k} \sum_{i=0}^{k} {k \brack i} (-1)^{k-i} T_{n+i-1}(x)
\end{equation}
Which is the inverse Stirling transform \cite{ref. 3} of following formula:
\begin{equation}
    T_{n+k-1}(x) = \sum_{i=0}^{k} {k \brace i} x^i Y_{n,i}(x)
\end{equation}
Using the definition of Touchard polynomial and (\ref{eq. 39}):
\begin{equation}
    T_{n+k-1}(x) = \sum_{i=0}^{k} \sum_{j=0}^{n} {k \brace i} {n+i \brace i+j}_i x^{i+j} = \sum_{m=0}^{n+k} (\sum_{i=0}^{k} {k \brace i}{n+i \brace m}_i) x^m
\end{equation}
We get:
\begin{equation}
    {n+k-1 \brace m} = \sum_{i=0}^{k} {k \brace i}{n+i \brace m}_i
\end{equation}

\subsection{Solving $y_{m,k}(x)=k y_{m-1,k}(x)+(x+k) y_{m-1,k+1}(x)$}

First, we consider $s_{n,k}(x) = \Gamma (x+k) y_{n,k}(x)  , s_{0,k}(x) = \Gamma (x+k)$  and use bottom-up unfolding as in 3.4:
\begin{equation}
    s_{n,k}(x) = \sum_{i=0}^{n} {n+k \brace i+k}_k \Gamma (x+k+i)
\end{equation}
Now by the definition of gamma function we get:
\begin{equation}
\begin{split}
    s_{n,k}(x) & = \sum_{i=0}^{n} {n+k \brace i+k}_k \int_{0}^{\infty} t^{x+k+i-1} e^{-t} dt \\
    & = \int_{0}^{\infty} \sum_{i=0}^{n} {n+k \brace i+k}_k t^{x+k+i-1} e^{-t} dt \\
    & = \int_{0}^{\infty} T_{n,k}(t) t^{x+k-1} e^{-t} dt
\end{split}
\end{equation}
Then, by $s_{n,k}(x) = \Gamma (x+k) y_{n,k}(x)$ we get:
\begin{equation}
    y_{n,k}(x) = \frac{1}{\Gamma (x+k)} \int_{0}^{\infty} T_{n,k}(t) t^{x+k-1} e^{-t} dt
\end{equation}
Finally, we have the explicit formula of the following rising factorial polynomial:
\begin{equation}
	\begin{split}
		\sum_{k=0}^{n} {n \brace k} x^{\overline{k}} & = y_{n,0}(x) \\
		& = \frac{1}{\Gamma (x)} \int_{0}^{\infty} T_n (t) t^{x-1} e^{-t} dt
	\end{split}
\end{equation}

\section{Rising and falling factorial transforms}

In this section, we discuss the transformation of power polynomials into rising and falling factorial polynomials. These transforms justify the results that we obtained in the previous section. First, we consider a rising factorial polynomial and a power polynomial with the same coefficient as follows respectively:
\begin{equation}
\label{eq. 55}
    R_n (x) = \sum_{k=0}^{n} a_{n,k} x^{\overline{k}}
\end{equation}
\begin{equation}
\label{eq. 56}
    P_n (x) = \sum_{k=0}^{n} a_{n,k} x^k
\end{equation}
We can apply the definition of rising factorial as below:
\begin{equation}
    R_n (x) = \sum_{k=0}^{n} a_{n,k} \frac{\Gamma (x+k)}{\Gamma (x)}
\end{equation}
As well as the definition of gamma function:
\begin{equation}
\begin{split}
\label{eq. 58}
    RFT(P_n (x)) & = R_n (x) \\
    & = \frac{1}{\Gamma (x)} \sum_{k=0}^{n} a_{n,k} \int_{0}^{\infty} t^{x+k-1} e^{-t} dt \\
    & = \frac{1}{\Gamma (x)} \int_{0}^{\infty} (\sum_{k=0}^{n} a_{n,k} t^k) t^{x-1} e^{-t} dt \\
    & = \frac{1}{\Gamma (x)} \int_{0}^{\infty} P_n (t) t^{x-1} e^{-t} dt
\end{split}
\end{equation}
Here we encounter the fact that operation can be done on two distinct types of polynomials (i.e. power and factorial polynomials). This is useful to derive formulas for $RFT(RFT(P_n (x)))$ and $RFT^{-1}(RFT^{-1}(P_n(x)))$ utilizing $x^n = \sum_{k=0}^{n} {n \brace k} (x)_k $ and $x^{\overline{n}} = \sum_{k=0}^{n} {n \brack k} x^k $ on (\ref{eq. 55}) and (\ref{eq. 56}), respectively:
\begin{equation}
\begin{split}
    P_n (-x) & = \sum_{k=0}^{n} a_{n,k} (-1)^k x^k = \sum_{k=0}^{n} \sum_{i=0}^{k} {k \brace i} a_{n,k} (-1)^k (x)_i \\
    & = \sum_{i=0}^{n} (-x)^{\overline{i}} \sum_{k=i}^{n} {k \brace i} a_{n,k} (-1)^{k-i}
\end{split}
\end{equation}
\begin{equation}
\begin{split}
    RFT^{-1} (RFT^{-1} (R_n (x))) & = RFT^{-1} (P_n (x)) \\
    & = \sum_{i=0}^{n} x^i \sum_{k=i}^{n} {k \brace i} a_{n,k} (-1)^{k-i} \\
    & = \sum_{k=0}^{n} a_{n,k} (-1)^k T_k (-x)
\end{split}
\end{equation}
\begin{equation}
    R_n (x) = \sum_{k=0}^{n} a_{n,k} x^{\overline{k}} = \sum_{k=0}^{n} \sum_{i=0}^{k} {k \brack i} a_{n,k} x^i = \sum_{i=0}^{n} x^i \sum_{k=i}^{n} {k \brack i} a_{n,k}
\end{equation}
\begin{equation}
\begin{split}
    RFT(RFT(P_n (x))) & = RFT(R_n (x)) = \sum_{i=0}^{n} x^{\overline{i}} \sum_{k=i}^{n} {k \brack i} a_{n,k} \\
    & = \sum_{k=0}^{n} a_{n,k} (\sum_{i=0}^{k} {k \brack i} x^{\overline{i}}) = \frac{1}{\Gamma (x)} \int_{0}^{\infty} R_n (t) t^{x-1} e^{-t} dt
\end{split}
\end{equation}
Consider the following polynomial to which this form $ \Delta^k f(x) = (E-I)^k f(x) = \sum_{i=0}^{k} \binom{k}{i} (-1)^{k-i} f(x+i)$ and Cauchy multiplication are applied:
\begin{equation}
\begin{split}
    \sum_{k=0}^{\infty} \frac{\Delta^k P_n (-x)}{k!} y^k & = \sum_{k=0}^{\infty} \frac{y^k}{k!} \sum_{i=0}^{k} \binom{k}{i} (-1)^{k-i} P_n (-x-i) \\
    & = \sum_{k=0}^{\infty} y^k \sum_{i=0}^{k} \frac{(-1)^{k-i}}{(k-i)!} \frac{P_n (-x-i)}{i!} \\
    & = e^{-y} \sum_{k=0}^{\infty} \frac{P_n (-x-k)}{k!} y^k
\end{split}
\end{equation}
This is generalization of Dobinski formula \cite{ref. 5}. Also, we can use $ \frac{\Delta^k f(x)}{k!} = \sum_{i=k}^{\infty} {i \brace k}  \frac{f^{(i)} (x)}{i!} $ from \cite{ref. 12} for above polynomial:
\begin{equation}
\begin{split}
    \sum_{k=0}^{\infty} \frac{\Delta^k P_n (-x)}{k!} y^k & = \sum_{k=0}^{\infty} y^k \sum_{i=k}^{\infty} {i \brace k} \frac{P_n^{(i)} (-x)}{i!} \\
    & = \sum_{i=0}^{\infty} \frac{P_n^{(i)} (-x)}{i!} \sum_{k=0}^{i} {i \brace k} y^k = \sum_{i=0}^{\infty} \frac{P_n^{(i)}(-x)}{i!} T_i (y)
\end{split}
\end{equation}
Consequently, by $ \frac{P_n^{(i)}(-x)}{i!} {\Bigg \rvert}_{x=0} = a_{n,k} (-1)^k $ we have:
\begin{equation}
    \sum_{i=0}^{n} a_{n,k} (-1)^k T_i (y) = \sum_{k=0}^{\infty} \frac{\Delta^k P_n (-x)}{k!} {\Bigg \rvert}_{x=0} y^k = e^{-y} \sum_{k=0}^{\infty} \frac{P_n (-k)}{k!} y^k
\end{equation}
We can derive inverse rising factorial transform:
\begin{equation}
\begin{split}
    RFT^{-1} (P_n (x)) & = \sum_{i=0}^{n} a_{n,k} (-1)^k T_i (-x) \\
    & = \sum_{k=0}^{\infty} \frac{\Delta^k P_n (-z)}{k!} {\Bigg \rvert}_{z=0} (-x)^k \\
    & = e^{x} \sum_{k=0}^{\infty} \frac{(-1)^k P_n (-k)}{k!} x^k
\end{split}
\end{equation}
A power and a falling factorial polynomial can be represented as follows:
\begin{equation}
    P_n (x) = \sum_{k=0}^{n} a_{n,k} x^k
\end{equation}
\begin{equation}
    F_n (x) = \sum_{k=0}^{n} a_{n,k} (x)_k
\end{equation}
\begin{equation}
    P_n (x) \xrightarrow{FFT}  F_n (x)
\end{equation}
By \cite{ref. 13}:
\begin{equation}
\label{eq. 70}
    \frac{1}{\Gamma (z)} = \frac{1}{2 \pi i} \int_{C} \frac{e^{t}}{t^{z}} dt
\end{equation}
Which $C$ is any deformed anticlockwise Bromwich contour \cite{ref. 13}. Using the definition of falling factorial:
\begin{equation}
    (x)_n = \frac{\Gamma (x+1)}{\Gamma (x-n+1)} = \frac{\Gamma (x+1)}{2 \pi i} \int_{C} \frac{e^{t}}{t^{x-n+1}} dt
\end{equation}
So, we get:
\begin{equation}
    F_n (x) = \frac{\Gamma (x+1)}{2 \pi i} \int_{C} \frac{P_n (t) e^{t}}{t^{x+1}} dt
\end{equation}
By \cite{ref. 14}:
\begin{equation}
    F_n (x) = FFT(P_n (x)) = \partial_{t}^{x} (P_n (t) e^{t}) {\Bigg \rvert}_{t=0}
\end{equation}
This formula demonstrate that falling factorial transform is equivalent to $x$th derivative of $P_n (t) e^{t}$.

(\ref{eq. 58}) is recall us fractional integral in $(0,\infty)$.
So we can rewrite (\ref{eq. 58}) for $RFT$ as follows:
\begin{equation}
    RFT(P_n (x)) = {}_{0}D^{-x}_{\infty} (P_n (t) e^{t})
\end{equation}
Also, There is another way to calculate $RFT$ owing to the definition of Mellin transform \cite{ref. 14}:
\begin{equation}
    RFT(P_n (x)) = \frac{1}{\Gamma (x)} \mathcal{M} (P_n (t) e^{-t})
\end{equation}

\section{Future works}

The following suggestions are made for further studies:
\begin{enumerate}
    \item Generalizing of this method: We proposed a method that was used to solve finite summation with a restrictive parameter. Hence, infinite summations with more variables can be considered. Also, this method can be generalized to multidimensional summation;
    \item Solving more examples: More examples can be solved by this method, especially polynomials (e.g., power or factorial polynomials);
    \item Extending rising and falling factorial transforms.
\end{enumerate}

\section{Conclusion}

We proposed a method that solves the summations that are finite and have a restrictive parameter. We illustrated some examples from which some useful identities emerge. We also investigated two new transforms that unify the results obtained in the example section. In addition, they are suitable for the transformation of the power series into the Newton series and vice versa. By the inversion formula of both transforms, almost all series can be solved. We can conclude from the results of section 4 that falling and rising factorial transforms can be represented as fractional derivatives and integrals of the intended function, respectively.

\end{document}